\documentclass[11pt]{amsart}
\usepackage{amssymb}
\usepackage{amsmath}
\usepackage{mathtools, float}
\usepackage{amsfonts}
\usepackage{stmaryrd}
\usepackage[english]{babel}
\usepackage[utf8]{inputenc}
\usepackage{enumerate}
\usepackage{graphicx,pinlabel,xcolor, fullpage}
\usepackage{adjustbox}
\usepackage{array,multirow}
\usepackage{float}

\theoremstyle{plain}
\newtheorem{theorem}{Theorem}[section]

\newtheorem{lemma}[theorem]{Lemma}

\newtheorem{theorem*}{Theorem}

\theoremstyle{definition}

\newtheorem{remark}[theorem]{Remark}
\newtheorem{note}[theorem]{Note}
\newtheorem{question}[theorem]{Question}
\newtheorem{proposition}[theorem]{Proposition}
\newtheoremstyle{case}{}{}{}{}{}{:}{ }{}
\theoremstyle{case}

\newcommand{\QED}{\hfill $\blacksquare$}

\begin{document}
\title[Uni-width subgroups, universal elements and lambda number of finite groups]{Uni-width subgroups, universal elements and lambda number of finite groups}
\author{Siddhartha Sarkar}
\address{Department of Mathematics\\
Indian Institute of Science Education and Research Bhopal\\
Bhopal Bypass Road, Bhauri \\
Bhopal 462 066, Madhya Pradesh\\
India}
\email{sidhu@iiserb.ac.in}



\begin{abstract}
A cyclic subgroup $N$ of a finite group $G$ is called a uni-width subgroup of $G$ if $N$ is the unique cyclic subgroup of $G$ of order $|N|$. In this article, we prove that a finite group $G$ admits a unique largest uni-width subgroup denoted by $U(1;G)$. We then show that the prime factors of the order of $U(1;G)$ influence the structure decomposition of its Fitting subgroup ${\mathrm{Fit}}(G)$. A power graph $\Gamma_G$ of a finite group is defined by $G$ being its set of vertices, and a pair of distinct elements $x,y \in G$ are connected by an edge if either $x \in \langle y \rangle$ or $y \in \langle x \rangle$. A universal element of a graph is a vertex that is adjacent to each of the remaining vertices. Our following result shows that a power graph $\Gamma_G$ of a finite non-trivial group admits a non-identity universal element if and only if it is either cyclic or a generalized quaternion $2$-group. The lambda number $\lambda(G)$ of a finite group $G$ is a measure of the least number of colors required for an $L(2,1)$-type of vertex coloring on $\Gamma_G$, which is known to be $\geq |G|$. Generalizing an earlier result, we then derive a necessary condition on a finite group $G$ such that $\lambda(G) = |G|$. Finally, we show that this result is best possible by exhibiting a family of groups without the necessary condition for which $\lambda(G) > |G|$.        
\end{abstract}


\maketitle

\section{Introduction}
\label{introsec}

\noindent All groups in this article are finite. For a prime $p \in {\mathbb{N}}$, by a $p$-group we mean a finite group $G$ with its order $|G| = p^n$ for some $n \in {\mathbb{N}}$. For a group $G$, we write $\Delta(G) := \{ |g| ~:~ g \in G \}$, where $|g|$ denote the order of the element $g \in G$. For $d \in \Delta(G)$, we define the {\it width} ${\mathfrak{m}}_{G}(d)$ of $d$ as the number of cyclic subgroups of $G$ of order $d$. An integer $d \in \Delta(G)$ is called {\it uni-width} if ${\mathfrak{m}}_{G}(d) = 1$ and in this case the associated cyclic subgroup, denoted by ${\mathcal{N}}_d \leq G$ of order $d$ is called an {\it uni-width subgroup} of $G$. 

\bigskip

\noindent The problem of counting ${\mathfrak{m}}_{G}(d)$ is classical in case $G$ is a finite $p$-group for some prime $p$ (see \cite{mil, kula}, \cite[$\S 5$]{ber} and \cite[Theorem 1.10, 1.17]{ber2} (for finite $2$-groups, such counting can be obtained as a consequence of a result of Alperin-Feit-Thompson, cf. \cite[Theorem 4.9]{isa}). Our interest in arbitrary non-trivial finite groups $G$ that admit non-trivial uni-width subgroups, relies on its connection to the question \ref{qn} on the power graph of $G$ stated below. 

\bigskip

\noindent An uni-width subgroup of $G$ is characteristic in $G$ and hence is contained in the Fitting subgroup ${\mathrm{Fit}}(G)$, which is the largest normal nilpotent subgroup of $G$. Recall that ${\mathrm{Fit}}(G) = \prod_{p \mid |G|, p {\mathrm{~prime}} } O_p(G)$ as a direct product, where $O_p(G)$ denotes the intersection of all Sylow-$p$-subgroups of $G$ for a prime $p \mid |G|$, called the $p$-core of $G$, and is a characteristic subgroup of $G$, since an automorphism of $G$ permutes its Sylow-$p$-subgroups. The existence of non-trivial uni-width subgroups influence the structure of ${\mathrm{Fit}}(G)$, which is our first result in this article: 

\bigskip

\begin{theorem}\label{uni-width-and-fitting}
Let $G$ be a non-trivial finite group, and let $1 \neq d \in \Delta(G)$ with ${\mathfrak{m}}_{G}(d) = 1$. Let $p \mid d$ be a prime factor. Then: 

\smallskip

\noindent (i) if $p$ is odd, then $O_p(G)$ is a cyclic subgroup of $G$. 

\smallskip

\noindent (ii) if $p=2$, then either ${\mathrm{O}}_2(G)$ is cyclic, or else ${\mathrm{O}}_2(G)$ is a $2$-group of maximal class and the Sylow-$2$-subgroup of ${\mathcal{N}}_d$ is contained in the (maximal) cyclic subgroup of index $2$ in ${\mathrm{O}}_2(G)$. 
\end{theorem}

\bigskip

\noindent As a consequence, we first prove the existence of the largest uni-width subgroup which contains all the uni-width subgroups of $G$: 

\bigskip

\begin{theorem}\label{uni-core-subgroup}
Let $U(1;G)$ be the subgroup of $G$ defined by 
\[
U(1;G) := \Big\langle x \in G ~:~ x \neq 1, ~{\mathfrak{m}}_{G}(|x|) = 1 \Big\rangle. 
\] 
Then, $U(1;G)$ is an uni-width subgroup, and consequently a characteristic subgroup of $G$. Moreover, for any uni-width subgroup ${\mathcal{N}}_d$ of $G$, we have ${\mathcal{N}}_d \leq U(1;G)$; i.e., $U(1;G)$ is the unique largest uni-width subgroup of $G$.   
\end{theorem}

\bigskip

\noindent We call $U(1;G) \leq G$ as the {\it uni-core} of $G$. Note that $U(1;G) = 1$ if $G$ admits no non-trivial uni-width subgroup. One interesting question at this point is as follows: Let $G$ be a non-trivial group and $\delta(G) := \{ p \in {\mathbb{N}} ~:~ p {\mathrm{~prime}}, O_p(G) {\mathrm{~is~cyclic}} \}$. Let us call $G$ to be {\it full uni-core} if either $U(1;G) = \prod_{p \in \delta(G)} O_p(G)$, or else $U(1;G)$ is an index $2$ subgroup of $\prod_{p \in \delta(G)} O_p(G) \cdot O_2(G)$. Note that, $G$ is full uni-core if and only if for each prime $p \mid \lvert U(1;G) \rvert$, the Sylow-$p$-subgroup $U_p(G)$ of $U(1;G)$ coincides with $O_p(G)$ for all odd $p$, and either $U_2(G) = O_2(G)$ when $O_2(G)$ is cyclic, or $[O_2(G) : U_2(G)] = 2$ when $O_2(G)$ is of maximal class.

\bigskip

\noindent Now using the methods of proof of Theorems \ref{uni-width-and-fitting} and \ref{uni-core-subgroup}, we can observe that a non-trivial nilpotent group is full uni-core. Also, if ${\mathrm{Fit}}(G) = 1$ (e.g., if $G$ is finite simple, or almost simple), then it is full uni-core. So, we can ask (see section \ref{final-example} for a family of examples):   

\bigskip

\begin{question}
Classify non-nilpotent full uni-core finite groups $G$ for which $U(1;G) \neq 1$.
\end{question}




\bigskip

\noindent {\bf Power graph and the universal elements of $G$:} Let $G$ be a non-trivial finite group. The (undirected) power graph $\Gamma_G$ of $G$ is the graph with the vertex set $V(\Gamma_G) := G$ and any two distinct elements $x, y \in G$ are connected by an edge if either $x \in \langle y \rangle$, or $y \in \langle x \rangle$. It can be observed that $\Gamma_G$ is a simple graph; i.e., it contains no loops and parallel edges. 

\bigskip

\noindent The idea of (directed) power graph of a semigroup $G$ was introduced by Kelarev and Quinn \cite{kq} and Chakrabarty ${\mathrm{\it{ et. al.}}}$ \cite{cgs} to resolve a classification problem. This article only considers the undirected power graphs corresponding to finite groups. Later on, in \cite{cam,cgh} it was shown that the power graph of two finite groups $G$ and $H$ with the same order, are isomorphic if and only if given any $d \mid |G|$, the number of elements of $G$ of order $d$ is equal to the number of elements of $H$ of order $d$. For a detailed survey on the power graphs and related areas, see the survey articles \cite{akc} and \cite{ksc}. 

\bigskip

\noindent An element $x \in G$ is called an {\it universal element} if $x$ is connected to every vertex in $G \setminus \{ x \}$ by an edge in $\Gamma_G$. Note that the identity element of $G$ is always a universal element. Our next result shows that the existence of a non-trivial universal element in a finite group is quite rare: 

\bigskip

\begin{theorem}\label{grps-with-univ-elms}
Let $G$ be a finite non-trivial group. Then, the following statements are equivalent:

\smallskip

\noindent (i) $G$ contains a non-trivial universal element,

\smallskip

\noindent (ii) $G$ is either cyclic or a generalized quaternion $2$-group.

\end{theorem}

\bigskip

\noindent {\bf Lambda number of a finite group:} Let $G$ be a finite group. An $L(2,1)$-labelling (or $L(2,1)$-coloring) of the graph $\Gamma_G$ is a function $f : G \rightarrow {\mathbb{Z}}$ such that $|f(x) - f(y)| \geq 2$ (resp. $\geq 1$) if $x$ and $y$ are adjacent (resp. if the shortest path between $x$ and $y$ in $\Gamma_G$ has length $2$).

\smallskip

\noindent For any $L(2,1)$-labelling of $\Gamma_G$, the span of $f$, denoted by ${\mathrm{span}}(f)$ is defined as
\[
{\mathrm{span}}(f) := \max_{x \in G} f(x) - \min_{x \in G} f(x).
\]
Let ${\mathcal{L}}(2,1)$ denote the set of all $L(2,1)$-labellings of $\Gamma_G$. The lambda number $\lambda(G)$ of a finite group $G$ is then defined by
\[
\lambda(G) := \min \Big\{ {\mathrm{span}}(f) ~;~ f \in {\mathcal{L}}(2,1) \Big\}.
\]

\bigskip

\noindent In \cite{mfw}, the following result was proved:

\begin{theorem}\label{mfw-main-lemma} (\cite[Lemma 3.1]{mfw})
Let $G$ be a finite group. Then, $\lambda(G) \geq |G|$ and the equality holds if and only if the complement graph $(\Gamma_G \setminus \{ 1 \})^c$ contains a Hamilton path; i.e., a linear subgraph that contains all the vertices of $G \setminus \{ 1 \}$. 
\end{theorem}

\bigskip

\noindent From this, the following question was posed:

\bigskip

\begin{question}\label{qn}
Classify all finite groups $G$ for which $\lambda(G) = |G|$.
\end{question}

\bigskip

\noindent So far, we know a few of the non-trivial finite groups $G$ with $\lambda(G) = |G|$ (see \cite{mfw}). In an earlier article, we proved the following result: 

\bigskip

\begin{theorem}\label{lambda-number-p-group} (\cite{ms2}) 
Let $G$ be a non-trivial finite group of prime power order. Then $\lambda(G) = |G|$ if and only if $G$ is neither cyclic, nor a generalized quaternion $2$-group. 
\end{theorem}

\bigskip

\noindent So far, the only known examples of non-trivial finite groups for which $\lambda(G) > |G|$ are either cyclic or a generalized quaternion. In light of these, it is a tempting question to find a necessary and sufficient criterion for finite groups over more general classes satisfying the equality $\lambda(G) = |G|$. The next result of this article shows that $\lambda(G) = |G|$ holds under some stronger hypothesis than the one of the equivalent condition in Theorem \ref{grps-with-univ-elms}.   

\bigskip

\noindent Let $G$ be a finite non-trivial group and let $d \in \Delta(G)$. Let $m := {\mathfrak{m}}_G(d)$ and let $C_1, \dotsc, C_m$ denotes the cyclic subgroups of $G$ of order $d$. We call $d$ {\it singular} if $C_i \cap U(1;G) \neq 1$ for some $1 \leq i \leq m$, and  $\lvert \{ j ~:~ C_j \cap U(1;G) = 1 \} \rvert \leq 2$.  The subset of $\Delta(G)$ consisting of all singular integers is denoted by $\Delta_0(G)$.

\bigskip

\begin{theorem}\label{main-thm-classification}
Let $G$ be a finite non-trivial group. 

\smallskip

\noindent (i) If $U(1;G) = 1$, then $\lambda(G) = |G|$. In particular, if $G$ is a finite non-abelian simple (resp., an almost simple) group, then $\lambda(G) = |G|$.

\smallskip

\noindent (ii) If $U(1;G) \neq 1$, let ${\mathfrak{S}} := \{ C \leq G ~:~ C {\mathrm{~cyclic}},~ |C| \in \Delta_0(G) \}$. Assume that
\[
\Big\lvert \bigcup_{C \in {\mathfrak{S}}} C \setminus \{ 1 \} \Big\rvert \leq \Big\lvert G \setminus \bigcup_{C \in {\mathfrak{S}}} C \Big\rvert. 
\]    
Then, $\lambda(G) = |G|$. In particular, if $\langle \bigcup_{C \in {\mathfrak{S}}} C \rangle$ is a proper subgroup of $G$, then $\lambda(G) = |G|$. 
\end{theorem}

\bigskip

\noindent Note that if $G$ either non-trivial cyclic, or a generalized quaternion $2$-group, then $U(1;G) = G$, and then $\bigcup_{C \in {\mathfrak{S}}} C = G$ (see Lemma \ref{uni-width-hull-normal}). This shows that the hypothesis of the above theorem \ref{main-thm-classification} is stronger than the hypothesis that $G$ is neither cyclic, nor a generalized quaternion $2$-group. 

\bigskip

\noindent In the end, we show that the hypothesis (ii) is optimum, and thereby construct an infinite family of non-trivial groups $G$ which is neither cyclic, nor a generalized quaternion, satisfying $\lambda(G) > |G|$, for which the hypothesis (ii) is false.

\bigskip

\noindent The article is organized as follows: In section \ref{sec-uni-width}, we discuss the structural properties of the uni-width subgroups and the proofs of Theorems \ref{uni-width-and-fitting} and \ref{uni-core-subgroup}. In section \ref{universal-elms}, we discuss the properties of universal elements and the proof of Theorem \ref{grps-with-univ-elms}. In section \ref{final-main-thm-sec} we present the proof of Theorem \ref{main-thm-classification}. In the final section \ref{final-example}, we provide an example of an infinite family of non-trivial finite groups satisfying $\lambda(G) > |G|$. 

\bigskip

\section{Uni-width subgroups}
\label{sec-uni-width}

\bigskip

\noindent For most of the part, we will write a decomposition of an arbitrary element $1 \neq x \in G$ into its prime power order components. So, it would be useful to write it down in the following form:

\bigskip

\begin{lemma}\label{elms-decomposition}
Let $G$ be a finite group and $1 \neq x \in G$ be arbitrary. Let $I_{x}$ denote the set of primes $p \mid |x|$. Then we can decompose $x = \prod_{p \in I_{x}} x_p$ with $x_p \in \langle x \rangle$ is of $p$-power order for every prime $p \in I_{x}$. If $H$ is a nilpotent subgroup of $G$ with $x \in H$, then this decomposition is unique in $H$. 
\end{lemma} 

\bigskip

\noindent {\bf Proof.} Without loss of generality, we may assume that $G$ is nilpotent. Let $J$ denotes the set of primes $q \mid |G|$, and we have $\{ p_1, \dotsc p_r \} := I_{x} \subseteq J$. If $|x| = d$, we can write $d = p^{\alpha_1}_1 \dotsc p^{\alpha_r}_r$. Define $d_i := p^{\alpha_1}_1 \dotsc {\widehat{p^{\alpha_i}_i}} \dotsc p^{\alpha_r}_r$ for every $1 \leq i \leq r$, where the symbol "$~\widehat{}~$" denotes the absence of the particular element. Then ${\mathrm{g.c.d}}(d_1, \dotsc, d_r) = 1$, and hence there exists $\lambda_1, \dotsc, \lambda_r \in {\mathbb{Z}}$ such that $\lambda_1 d_1 + \dotsc + \lambda_r d_r = 1$ and $\lambda_i \not\equiv 0$ mod $p_i$ for every $1 \leq i \leq r$. Now we have
\[
x = x^{\lambda_1 d_1} \dotsc x^{\lambda_r d_r}
\]
and from above conditions, we have $|x^{\lambda_i d_i}| = p^{\alpha_i}_i$ for every $1 \leq i \leq r$. Setting $x_{p_i} := x^{\lambda_i d_i} ~(1 \leq i \leq r)$, we obtain the decomposition. Now, since $G$ is nilpotent, it is a direct product of its unique Sylow-subgroups, and the component $x_{p_i}$ belongs to the Sylow-$p_i$-subgroup, from which the uniqueness follows. \QED 

\bigskip

\begin{note}\label{projection-comps}
In case $G$ is a finite non-trivial nilpotent group and $1 \neq x \in G$, we may consider the projections $\pi_p : G \rightarrow G_p$ for each prime $p \in J$ in the notations of above lemma \ref{elms-decomposition}, and define $\pi_p(x) := x_p$. Using the unique Sylow-subgroup decomposition of $G$, it follows that these maps are well-defined surjective homomorphisms, with $\pi_p(x) = 1$ if and only if $p \nmid |x|$.  
\end{note}

\bigskip

\begin{note}\label{p-grps-max-class}
Before we prove Theorem \ref{uni-width-and-fitting}, we need to recall the following definition: Let $p \in {\mathbb{N}}$ be a prime and $G$ be a finite group of order $p^n$ for some $n \in {\mathbb{N}}$. It is well known that $G$ is nilpotent with nilpotency class $\leq n-1$. Then, $G$ is said to be a $p$-group of maximal class if its nilpotency class is equal to $n-1$. If $p=2$, then the $2$-groups of maximal class fall into three mutually disjoint families (see \cite[Theorem 5.1]{falc}):

\smallskip

\noindent {\bf Dihedral groups.}
\[
{\mathbb{D}}_{2^{e+1}} = \Big\langle a, b ~:~ a^{2^e} = 1, b^2 = 1, b^{-1} ab = a^{-1} \Big\rangle \hspace{.5in} (e \geq 2),
\]
\noindent {\bf Generalized Quaternion groups.}
\[
{\mathbb{Q}}_{2^{e+1}} = \Big\langle a, b ~:~ a^{2^e} = 1, b^2 = a^{2^{e-1}}, b^{-1} ab = a^{-1} \Big\rangle \hspace{.5in} (e \geq 2),
\]
\noindent {\bf Semi Dihedral groups.}
\[
{\mathbb{SD}}_{2^{e+1}} = \Big\langle a, b ~:~ a^{2^e} = 1, b^2 = 1, b^{-1} ab = a^{2^{e-1} -1} \Big\rangle \hspace{.5in} (e \geq 3).
\]
These three families of $2$-groups of maximal class always contain a unique maximal cyclic subgroup of index $2$, namely $\langle a \rangle$.  
\end{note}

\bigskip

\noindent {\bf Proof of Theorem \ref{uni-width-and-fitting}.} (i) Let ${\mathcal{N}}_d = \langle x \rangle$ and let $p \mid d$ be an odd prime. We decompose $x = x^{\prime} x^{\prime \prime}$ for some $x^{\prime}, x^{\prime \prime} \in {\mathcal{N}}_d$ so that $x^{\prime}$ is of $p$-power order and ${\mathrm{g.c.d.}}(|x^{\prime \prime}|, p) = 1$. From (ii), we have $x^{\prime} \in {\mathrm{O}}_p(G)$. If ${\mathrm{O}}_p(G)$ is not cyclic, then using Miller's theorem \cite[Theorem 1.10(b)]{ber2}, there exists $y \in {\mathrm{O}}_p(G)$ such that $|y| = |x^{\prime}|$ and $\langle y \rangle \neq \langle x^{\prime} \rangle$. Since $\langle y \rangle$ and $\langle x^{\prime \prime} \rangle$ are of co-prime order and $y$ centralizes $\langle x^{\prime \prime} \rangle$, since ${\mathrm{O}}_p(G)$ centralizes the remaining Sylow subgroups of ${\mathrm{Fit}}(G)$. Then $yx^{\prime \prime}$ has order $d$. From hypothesis, we have $yx^{\prime \prime} \in \langle x^{\prime} x^{\prime \prime} \rangle$, and consequently, $\langle y \rangle = \langle x^{\prime \prime} \rangle$, a contradiction.

\smallskip

\noindent (ii) Now we consider $p=2$. In case ${\mathrm{O}}_2(G)$ is neither cyclic, nor a $2$-group of maximal class, the proof is the same as above using \cite[Theorem 1.17]{ber2}. So we assume that ${\mathrm{O}}_2(G)$ is a $2$-group of maximal class with index $2$ normal cyclic subgroup $H \leq {\mathrm{O}}_2(G)$. If $x^{\prime} \not\in H$, then using \cite[Corollary 2.4]{ms2}, there exists $y \not\in H$ such that $\langle y \rangle \neq \langle x^{\prime} \rangle$, and $|y| = |x^{\prime}|$. An identical argument as above, as in the case of an odd prime, proves that $\langle y \rangle = \langle x^{\prime \prime} \rangle$, a contradiction. Hence we must have $x^{\prime} \in H$. \QED

\bigskip

\noindent The proof of Theorem \ref{uni-core-subgroup} require the following lemma:

\bigskip

\begin{lemma}\label{prod-uni-width-subgrps}
Let $G$ be a non-trivial finite group. Then, the product of two uni-width subgroups of $G$ is a uni-width subgroup of $G$.  
\end{lemma}

\bigskip

\noindent {\bf Proof.} Let $E = {\mathcal{N}}_{d_1}$ and $F = {\mathcal{N}}_{d_2}$ be two uni-width subgroups of $G$. Let $I_{E}$ (resp. $I_{F}$) denote the set of all odd prime factors of $d_1$ (resp. $d_2$) and we write $I := I_{E} \cup I_{F}$. Since $E, F, EF \leq {\mathrm{Fit}}(G)$, there are direct product decompositions 
\[
E = E_2 \cdot \prod_{p \in I_{E}} E_p, \hspace*{.2in} F = F_2 \cdot \prod_{p \in I_{F}} F_p, 
\]
and combining these, we have 
\[
EF = (E_2 F_2) \cdot \prod_{p \in I} E_p F_p.
\]
Now from Theorem \ref{uni-width-and-fitting}, for each prime $p \in I$ (resp. for $p=2$), we have $E_p$ and $F_p$ (resp. $E_2$ and $F_2$) are contained in the cyclic group ${\mathrm{O}}_p(G)$ (resp. in a cyclic subgroup of ${\mathrm{O}}_2(G)$) of prime power order. This shows that $E_p F_p = E_p \cup F_p$ is the (unique) cyclic Sylow-$p$-subgroup of $EF$ for every $p \in I \cup \{ 2 \}$, and the above equation is a direct product decomposition. Consequently, $EF$ is a cyclic group of order $d := {\mathrm{l.c.m.}}(d_1, d_2)$.

\smallskip

\noindent It remains to prove that $EF$ is a uni-width subgroup of $G$. Let $g \in G$ with $|g| = d$. We decompose $g = \prod_{q \in I \cup \{ 2 \}} g_q$ with elements $g_q \in \langle g \rangle$ of $q$-power order. For each $q \in I \cup \{ 2 \}$, we define a pair of elements $(a_q, b_q)$ as follows:
\[
(a_q, b_q) =
\begin{cases}
\big( g_q, g^{|E_q|/|F_q|}_q \big) &\mbox{if } F_q \subseteq E_q, \\ 
\big( g^{|F_q|/|E_q|}_q, g_q \big)	&\mbox{if } E_q \subset F_q.
\end{cases}
\] 
Now set $a := \prod_{q \in I \cup \{ 2 \}} a_q$ and $b := \prod_{q \in I \cup \{ 2 \}} b_q$. Then, $|a| = d_1, |b| = d_2$, and hence $\langle a \rangle = E, \langle b \rangle = F$. Using the similar arguments as above, it follows that $\langle a, b \rangle = EF$. Now, since $\langle a, b \rangle \subseteq \langle g \rangle$, and both of them have order $d$, the statement follows. \QED  

\bigskip

\begin{remark}\label{dihedral}
Note that a subgroup of a uni-width subgroup need not be a uni-width subgroup. For example, consider the dihedral group $G := {\mathbb{D}}_{2^{e+1}}$ of order $2^{e+1}$ as given in Note \ref{p-grps-max-class} for some integer $e \geq 2$. Then $\langle a \rangle \leq G$ is a uni-width subgroup of order $2^e$. However, the subgroup $Z(G) = \langle a^{2^{e-1}} \rangle$ of $\langle a \rangle$ of order $2$ is not uni-width, since ${\mathfrak{m}}_G(2) = 1 + 2^e$, which can be observed by counting the complete list of order $2$ subgroups: $\langle a^{2^{e-1}} \rangle, \langle ba^i \rangle ~(0 \leq i \leq 2^e - 1)$.  
\end{remark} 

\bigskip

\noindent {\bf Proof of Theorem \ref{uni-core-subgroup}.} If $U(1;G) = 1$, then by definition ${\mathfrak{m}}(|g|) \geq 2$ for all $g \in G \setminus \{ 1 \}$, and statement follows. So, assume that $U(1;G) \neq 1$, and let $X := \{ x \in G ~:~ x \neq 1, ~{\mathfrak{m}}_{G}(|x|) = 1 \} = \{ x_1, x_2, \dotsc, x_r \}$. Then, $U(1;G) = \langle X \rangle$. Next, by using the Lemma \ref{prod-uni-width-subgrps} repeatedly, we have $\langle X \rangle = \langle x_1 \rangle \langle x_2 \rangle \dotsc \langle x_r \rangle$ is an uni-width subgroup of $G$. Finally, for any non-trivial uni-width subgroup $\langle y \rangle = {\mathcal{N}}_d \leq G$, we have $y \in X$, and hence ${\mathcal{N}}_d \subseteq U(1;G)$. \QED

\bigskip

\section{Universal elements in finite groups}
\label{universal-elms}

\bigskip

\noindent In this section, we will prove Theorem \ref{grps-with-univ-elms}. Recall that, from the definition of a power graph of a finite group $G \neq 1$, we have $x \in G$ is an universal element of $G$, if for any $g \in G$, either $g = x^n$, or $x = g^m$ for some $m,n \in {\mathbb{N}}$. We begin with the following easy Lemma:

\bigskip

\begin{lemma}\label{univ-elms-central}
Let $G$ be a finite group and $x \in G$ be a universal element of $G$. Then $x \in Z(G)$, and consequently, $x$ belongs to the Fitting subgroup ${\mathrm{Fit}}(G)$ of $G$.
\end{lemma}

\bigskip

\noindent {\bf Proof.} Let $1 \neq g \in G$. If $g = x^n$, or $x = g^m$ for some integer $m,n \in {\mathbb{N}}$, then $[g,x] = 1$. Hence $x$ centralizes every element of $G$. Since $Z(G) \subseteq {\mathrm{Fit}}(G)$, the statement follows. \QED  

\bigskip

\noindent Next, we prove a stronger version of Theorem \ref{grps-with-univ-elms}:

\bigskip

\begin{proposition}\label{p-grps-with-univ-elms}
Let $p$ be a prime, and $G$ be a non-trivial finite $p$-group. Then $G$ contains a non-trivial universal element if and only if $G$ is either cyclic, or the generalized quaternion $2$-group. 
\end{proposition}

\bigskip

\noindent {\bf Proof.} "$\Leftarrow$" If $G$ is cyclic, then every generator of $G$ is an universal element. In case $G = {\mathbb{Q}}_{2^{e+1}}$ of order $2^{e+1}$ for some $e \geq 2$, then $Z(G) = \langle b^2 \rangle$ is cyclic of order $2$. Now, $b^2$ is the only element of $G$ of order $2$, and hence is a power of every non-trivial element of $G$. This proves that $b^2$ is a non-trivial universal element of $G$.

\smallskip

\noindent "$\Rightarrow$" So we assume that $G$ is a finite $p$-group and $1 \neq u \in G$ is an universal element of $G$. Now if $g \in G$, and $|g| = |u|$, then using the definition of universal element, we have $\langle g \rangle = \langle u \rangle$. This implies that ${\mathfrak{m}}(|u|) = 1$. If $p$ is odd, using \cite[Theorem 1.10]{ber2} it follows that $G$ is cyclic. Now suppose $p=2$, and $G$ is not cyclic. Then, from \cite[Theorem 1.17]{ber2}, this implies that $G$ must be a $2$-group of maximal class, and hence $G$ must be one of the groups ${\mathbb{D}}_{2^{e+1}}$, ${\mathbb{Q}}_{2^{e+1}}$ for some $e \geq 2$, or else ${\mathbb{SD}}_{2^{e+1}}$ for some $e \geq 3$ as given in Note \ref{p-grps-max-class}. 

\smallskip

\noindent First suppose $G = {\mathbb{D}}_{2^{e+1}}$. Since every element outside the maximal cyclic subgroup $\langle a \rangle$ has order $2$, we have $|u| \geq 4$ and $u \in \langle a \rangle$. Now consider the element $b \in {\mathbb{D}}_{2^{e+1}} \setminus \langle a \rangle$ of order $2$. Since $u$ is universal, there exists $m \in {\mathbb{N}}$ such that $u^m = b$, a contradiction. Hence, $G$ cannot be dihedral.

\smallskip

\noindent In case $G = {\mathbb{SD}}_{2^{e+1}}$, a similar argument shows that $|u| \geq 8$ and it belongs to the index $2$ cyclic subgroup of $\langle a \rangle \leq {\mathbb{SD}}_{2^{e+1}}$. However, any element outside $\langle a \rangle$ has order $2$, or $4$, and cannot be written as a power of $u$. Hence $G$ must be the generalized quaternion $2$-group. \QED    

\bigskip

\noindent {\bf Proof of Theorem \ref{grps-with-univ-elms}.} The part "$(ii) \Rightarrow (i)$" is done in the proof of Proposition \ref{p-grps-with-univ-elms}. So, we prove "$(i) \Rightarrow (ii)$". We assume that $G$ is a finite group with a universal element $u \neq 1$. 

\smallskip

\noindent We claim that $G$ is nilpotent. Indeed, if $G$ is not nilpotent, then ${\mathrm{Fit}}(G)$ is a proper subgroup of $G$. This implies that there exists a prime $p \mid |G|$, such that ${\mathrm{O}}_p(G)$ is a proper subgroup of the Sylow-$p$-subgroups of $G$. Let $H$ denote a Sylow-$p$-subgroup of $G$, and consider any $g \in H \setminus {\mathrm{O}}_p(G)$. Now, from Lemma \ref{univ-elms-central} we have $u \in {\mathrm{Fit}}(G)$ and $g \not\in {\mathrm{Fit}}(G)$. From definition of universal element, there exists $m \in {\mathbb{N}}$, such that $g^m = u$. This implies that $u$ is of $p$-power order. Since $G$ is not nilpotent, there exists a prime $q \neq p$ such that $q \mid |G|$. Now let $K$ be a Sylow-$q$-subgroup of $G$ and consider any $1 \neq h \in K$. Again, using the definition of universal element, either $h^n = u$, or $h = u^m$ for some integer $m,n \in {\mathbb{N}}$. This is a contradiction, since $h$ and $u$ have coprime orders. This proves that $G$ is nilpotent. 

\smallskip

\noindent Now we can write $G = \prod_{q \in I} G_q$ as a direct product of Sylow-subgroups, where $I$ is the set of primes $q \mid |G|$. Using the notations as in note \ref{projection-comps}, we decompose $u = \prod_{q\in I} \pi_q(u)$, and we claim that $\pi_p(u) \neq 1$ for each $p \in I$. If $\pi_p(u) = 1$ for some $p \in I$, then ${\mathrm{g.c.d.}}(|u|, p) = 1$. Now if $1 \neq g \in G_p$, then from the definition of universal element we have either $g^n = u$, or $g = u^m$ for some $m,n \in {\mathbb{N}}$. Since $u$ and $g$ are non-trivial, this is a contradiction. This proves that the projection components of $u$ are non-trivial in each of the Sylow subgroups of $G$.

\smallskip

\noindent Next we claim that $\pi_p(u)$ is an universal element in $G_p$ for every $p \in I$. To see this, fix any $p \in I$ and let $h \in G_p$. Then we have either $h^n = u$, or $h = u^m$ for some $m,n \in {\mathbb{N}}$. Applying the projection to each of these equations, we have either $h^n = \pi_p(u)$, or $h = \pi_p(u)^m$ for some $m,n \in {\mathbb{N}}$. Using the conclusion of the previous paragraph along with this, it shows that $\pi_p(u) \neq 1$ is a universal element in $G_p$ for every $p \in I$. 

\smallskip

\noindent Now suppose $G$ is not cyclic. Then using Proposition \ref{p-grps-with-univ-elms}, we have $2 \in I$, $G_2 = {\mathbb{Q}}_{2^{e+1}}$ for some $e \geq 2$, and $G_p$ is cyclic for every $p \in I \setminus \{ 2 \}$. It only remains to prove that $I = \{ 2 \}$. Assuming this is not the case, we fix the notations $I^{\prime} := I \setminus \{ 2 \}$,  
\[
G_2 = \Big\langle x_2, y_2 ~:~ x^{2^e}_2 = 1, x^{2^{e-1}}_2 = y^2_2, y^{-1}_2 x_2 y_2 = x^{-1}_2 \Big\rangle
\]
and $G_p = \langle x_p ~:~ x^{p^{e_p}}_p = 1 \rangle$ for some $e_p \geq 1$, and for each $p \in I \setminus \{ 2 \}$. Since the only universal element in $G_2$ is the unique order $2$ element, we have $\pi_2(u) = x^{2^{e-1}}_2 = y^2_2$. 

\smallskip

\noindent We claim that $|\pi_p(u)| = p$ for all $p \in I^{\prime}$. Consider the element $g = x_2 \prod_{p \in I^{\prime}} x^{p^{e_p - 1}}_p$. If $g = u^m$ for some $m \in {\mathbb{N}}$, then applying the projection $\pi_2$ we obtain $x_2 = \pi_2(u)^m$, a contradiction. Hence we must have $g^n = u$ for some $n \in {\mathbb{N}}$. Again applying the projection $\pi_p$ to this relation we have $\big( x^{p^{e_p - 1}}_p \big)^n = \pi_p(u)$. From this we must have ${\mathrm{g.c.d.}}(n,p) = 1$ and therefore, $|\pi_p(u)| = p$ for each $p \in I^{\prime}$. This shows that $|u| = 2 \prod_{p \in I^{\prime}} p$.

\smallskip

\noindent Now consider the element $x_2 \in G$. Then $x_2 \not\in \langle u \rangle$ and $u \not\in \langle x_2 \rangle$. This contradicts that $u$ is a universal element of $G$. Hence we have $I = \{ 2 \}$, which proves the theorem. \QED

\bigskip

\section{Cyclic classes and proof of Theorem \ref{main-thm-classification}}
\label{final-main-thm-sec}

\bigskip

\noindent We begin with a finer version of counting width as mentioned in section \ref{introsec}. Let $G$ be a non-trivial finite group and $d \in \Delta(G)$. For any $x,y \in G$, we define $x \sim y$ if $\langle x \rangle = \langle y \rangle$. It is easy to check that this is an equivalence relation defined on $G$. An equivalence class of this relation, called a {\it cyclic class}, consists of the elements which are the generators of a fixed cyclic subgroup of $G$. For any $d \in \Delta(G)$, we define
\[  
\Lambda(d;G) := \Big\{ g \in G ~:~ |g| = d \Big\}.
\]
Then $\Lambda(d;G)$ is equal to the disjoint union of all cyclic classes $[g]$ (represented by an element $g \in G$) so that $|g| = d$. For $g \in \Lambda(d;G)$, the number of elements of $G$ that belong to the cyclic class $[g]$ is $\phi(d)$, where $\phi$ denotes Euler's phi function. Then, we have ${\mathfrak{m}}_{G}(d) = \lvert \Lambda(d;G) \rvert /{\phi(d)}$. 

\bigskip

\noindent We leave a simple observation, which is used in multiple places in the rest of the article:

\begin{note}
Let $G$ be a non-trivial finite group, and $x,y \in G \setminus \{ 1 \}$. If $|x| \nmid |y|$ and $|y| \nmid |x|$, then $x$ is adjacent to $y$ in the complement graph $(\Gamma_{G} \setminus \{ 1 \})^c$. The same conclusion holds true if $|x| = |y|$ and $\langle x \rangle \neq \langle y \rangle$.  
\end{note} 

\bigskip

\noindent Let $d \in \Delta(G)$ with $m := {\mathfrak{m}}_{G}(d)$, and let $F_1, \dotsc, F_m$ denote the cyclic classes of $G$ consisting of elements of order $d$. Then, we have $\Lambda(d;G) = F_1 \cup \dotsc \cup F_m$. For any $1 \leq i, j \leq m$ with $i \neq j$, if $x \in F_i$ and $y \in F_j$ are arbitrarily chosen, then $x$ is adjacent to $y$ in the complement graph $(\Gamma_{G} \setminus \{ 1 \})^c$. Also, each $F_i$ form a discrete subgraph of $(\Gamma_{G} \setminus \{ 1 \})^c$. This shows that the set of vertices $F_1, \dotsc, F_m$ form an $m$-partition and induce an $m$-partite subgraph of $(\Gamma_{G} \setminus \{ 1 \})^c$.  

\bigskip

\noindent Now consider any sub-collection of $\{ F_1, \dotsc, F_m \}$ with at least two of $F_i$'s, which we can denote by $\{ F_1, \dotsc, F_k \}$ without any loss of generality. Now writing $F_i := \{ v_{i1}, v_{i2}, \dotsc, v_{i \phi(d)} \}$, we define $\gamma(F_1, \dotsc, F_k)$ to be the path 
\[
\gamma(F_1, \dotsc, F_k) := \big( v_{11}, v_{21}, \dotsc, v_{k1}, v_{12}, v_{22}, \dotsc, v_{k2}, \dotsc, v_{1 \phi(d)}, v_{2 \phi(d)}, \dotsc, v_{k \phi(d)} \big)
\]
in the complement graph $(\Gamma_{G} \setminus \{ 1 \})^c$. Notice that the initial vertex $v_{11}$ (resp. the terminal vertex $v_{k \phi(d)}$) can be chosen arbitrarily from $F_1$ (resp. from $F_k$) and the parts $F_1, \dotsc, F_k$ can also be arbitrarily ordered. This implies that, for any arbitrary $x,y \in F_1 \cup \dotsc \cup F_k$ that come from disjoint cyclic classes, there is a Hamiltonian path in the induced subgraph with vertices $F_1 \cup \dotsc \cup F_k$ that begins at $x$ and ends at $y$.     

\bigskip

\noindent Now let $d \in \Delta(G)$ and let $C_1, \dotsc, C_m$ denotes the cyclic subgroups of $G$ of order $d$, where $m = {\mathfrak{m}}_G(d)$. Recall that, $d$ is called singular if $U(1;G) \cap C_i \neq 1$ for some $1 \leq i \leq m$, and $\lvert \{ j ~:~ 1 \leq j \leq m, C_j \cap U(1;G) = 1 \} \rvert \leq 2$. The subset  $\Delta_0(G) \subseteq \Delta(G)$ is the set of singular elements. Recall that $S \subseteq G$ is called a characteristic subset if it is invariant under all automorphisms of $G$. A subgroup generated by a characteristic subset is a characteristic subgroup.   

\bigskip

\begin{lemma}\label{uni-width-hull-normal}
Let $G$ be a finite non-trivial group with $\Delta_0(G) \neq \emptyset$. Let ${\mathfrak{S}}$ denote the collection of all cyclic subgroups $C$ of $G$ such that $|C| \in \Delta_0(G)$. Then, $\cup_{C \in {\mathfrak{S}}} C$ is a characteristic subset of $G$ containing $U(1;G)$. Consequently, $\langle \cup_{C \in {\mathfrak{S}}} C \rangle$ is a characteristic subgroup of $G$ containing $U(1;G)$.
\end{lemma}

\bigskip

\noindent {\bf Proof.} Since $\Delta_0(G) \neq \emptyset$, we have $U(1;G) \neq 1$ in $G$. Now, by definition $U(1;G) \in {\mathfrak{S}}$ and hence $U(1;G) \subseteq \cup_{C \in {\mathfrak{S}}} C$. Next, for any $d \in \Delta(G)$, we have $\cup_{C ~{\mathrm{cyclic}}, |C| = d} C$ is a characteristic subset of $G$. Consequently, it follows that the union $\cup_{C \in {\mathfrak{S}}} C = \cup_{d \in \Delta_0(G)} \cup_{C ~{\mathrm{cyclic}}, |C| = d} C$ is a characteristic subset of $G$. \QED

\bigskip

\noindent {\bf Proof of Theorem \ref{main-thm-classification}.} In either of (i) or (ii), it is enough to construct a Hamiltonian path in $(\Gamma_{G} \setminus \{ 1 \})^c$.

\bigskip

\noindent {\bf Step 1.}

\bigskip

\noindent We begin with the observation that, if $G$ is a non-trivial group with $U(1;G) = 1$, then $\Delta_0(G) = \emptyset$. In this case, we set ${\mathfrak{S}} := \emptyset$ and $\cup_{C \in {\mathfrak{S}}} C = \{ 1 \}$. Now, we first construct a Hamilton path in $(\Gamma_{G} \setminus \{ 1 \})^c$ involving the vertices of $G \setminus \cup_{C \in {\mathfrak{S}}} C$. Note that for any $g \in G \setminus \cup_{C \in {\mathfrak{S}}} C$, we have $|g| \in \Delta(G) \setminus \Delta_0(G)$. 

\bigskip

\noindent For $d \in \Delta(G) \setminus \Delta_0(G)$, let $F_1, \dotsc, F_m$ denote the cyclic classes of order $d$ such that $\cup^m_{i=1} F_i \subseteq G \setminus \cup_{C \in {\mathfrak{S}}} C$. Since $d$ is not singular, $m \geq 3$. As noticed earlier at the beginning of this section, the induced subgraph of $(\Gamma_{G} \setminus \{ 1 \})^c$ consisting of the vertex set $F_1 \cup \dotsc \cup F_m$ contains a Hamiltonian path $\gamma_d := \gamma(F_1, \dotsc, F_m)$. The initial and terminal vertices of this path can be arbitrarily chosen from two distinct parts $F_1, \dotsc, F_m$. 

\bigskip

\noindent We write $\Delta(G) \setminus \Delta_0(G) = \{ d_1, d_2, \dotsc, d_k \}$ and we assume that $d_1 < d_2 < \dotsc < d_k$. For any $1 \leq i \leq k-1$, if $d_i \nmid d_{i+1}$, then the final vertex of $\gamma_{d_i}$ can be joined to the initial vertex of $\gamma_{d_{i+1}}$ by an edge in the complement graph $(\Gamma_{G} \setminus \{ 1 \})^c$. So assume that $d_i \mid d_{i+1}$. In this case, let $w$ denote the initial vertex of $\gamma_{d_{i+1}}$. Since a cyclic group of order $d_{i+1}$ has an unique subgroup of $d_i$ and the number of cyclic classes from $G \setminus \cup_{C \in {\mathfrak{S}}} C$ consisting of elements of order $d_i$ is $\geq 2$ ($\geq 3$ in case $U(1;G) \neq 1$), we may interchange the ordering of the cyclic classes in the definition of $\gamma_{d_i}$, to obtain the final vertex which is connected to $w$ in $(\Gamma_{G} \setminus \{ 1 \})^c$. Now, using the backward induction on $i < k$, we can construct the Hamilton path ${\mathcal{H}} := \gamma_{d_1} \gamma_{d_2} \dotsc \gamma_{d_k}$ which cover every vertex of $G \setminus \cup_{C \in {\mathfrak{S}}} C$. Notice that the number of cyclic classes involved in $\gamma_{d_1}$ is $\geq 2$ ($\geq 3$ in case $U(1;G) \neq 1$), among which only one class is adjacent to $\gamma_{d_2}$ in the Hamiltonian path ${\mathcal{H}}$. Thus, when we extend this to a Hamilton path in $(\Gamma_{G} \setminus \{ 1 \})^c$, if necessary, we can permute the remaining cyclic classes in $\gamma_{d_1}$ to obtain a suitable initial vertex of ${\mathcal{H}}$.

\bigskip

\noindent Now we can prove (i): If $U(1;G) = 1$, then the Hamiltonian path ${\mathcal{H}}^{\prime}$ covers all vertices of $(\Gamma_{G} \setminus \{ 1 \})^c$. This implies that, $\lambda(G) = |G|$. In case $G$ is a non-abelian simple group, we have $U(1;G) = 1$ from Theorem \ref{uni-core-subgroup}. Next, assume that $G$ is an almost simple group. Then we can write $S \leq G \leq {\mathrm{Aut}}(S)$ for some non-abelian finite simple group $S$, where in the inclusion we have the identification $S \equiv {\mathrm{Inn}}(S)$. It is also well known that $S = {\mathrm{Soc}}(G)$, where the socle ${\mathrm{Soc}}(G)$ of $G$ is defined to be the subgroup generated by all the minimal normal subgroups. We now prove that $U(1;G)$ is trivial: Assume, if possible $U(1;G) \neq 1$. Then, ${\mathrm{Fit}}(G) \neq 1$ and hence there exists a prime $p \mid |G|$ such that $O_p(G) \neq 1$. Since $O_p(G)$ is a finite $p$-group, it admits a non-trivial center. Now, from the following characteristic subgroup tower
\[
Z(O_p(G)) ~{\mathrm{char}}~ O_p(G) ~{\mathrm{char}}~ {\mathrm{Fit}}(G) ~{\mathrm{char}}~ G 
\] 
it follows that $G$ contains a non-trivial abelian minimal normal subgroup, which is inside the finite simple subgroup $S$ of $G$, by definition of socle. This is a contradiction, and hence $U(1;G) = 1$.

\bigskip

\noindent To establish part (ii), we now extend the path ${\mathcal{H}}$ to construct a Hamiltonian path of $(\Gamma_{G} \setminus \{ 1 \})^c$.
  
\bigskip

\noindent {\bf Step 2.}

\bigskip

\noindent Now we consider the vertices of $\cup_{C \in {\mathfrak{S}}} C$. We extend the Hamiltonian path ${\mathcal{H}}$ on its left end using the vertices in $\cup_{C \in {\mathfrak{S}}} C$. We first consider the vertices $g \in \cup_{C \in {\mathfrak{S}}} C$ that satisfy the following conditions:

\bigskip

\noindent (a) ${\mathfrak{m}}_G(|g|) \geq 2$, \\
\noindent (b) number of cyclic classes consisting of elements of order $|g|$ inside $\cup_{C \in {\mathfrak{S}}} C$ is $\geq 2$ 

\bigskip

\noindent We assume that $e_1 < e_2 < \dotsc < e_m$ be the orders of these elements. For each $e \in \{ e_1, \dotsc, e_m \}$, if $W_1, \dotsc, W_{k_e}$ are the cyclic classes of elements of order $e$ inside $\cup_{C \in {\mathfrak{S}}} C$, then we can construct a path $\delta_e := \gamma(W_1, \dotsc, W_{k_e})$, as before. Now, if $e_m \nmid d_1$ and $d_1 \nmid e_m$, then we can join $\delta_{e_m}$ with $\gamma_{d_1}$. Also, if $e_m = d_1$, then the Hamilton path ${\mathcal{H}}$ can be easily extended to the remaining vertices of order $d_1$ in $\cup_{C \in {\mathfrak{S}}} C$ and thereby constructing ${\mathcal{H}}^{\prime} := \delta_{e_1} \dotsc \delta_{e_m} {\mathcal{H}}$ . So we assume $e_m \neq d_1$, and consider the following two cases:

\bigskip

\noindent {\bf Case I.} $e_m \mid d_1$

\bigskip

\noindent Let $J_1, J_2 \subseteq \cup_{C \in {\mathfrak{S}}} C$ denote two cyclic classes consists of elements of order $e_m$; $g_1 \in J_1$ and $g_2 \in J_2$. Let $w$ be the initial vertex of $\gamma_{d_1}$ so that $|w| = d_1$. Now, if $g_1 \not\in \langle w \rangle$, then $g_1$ is adjacent to $w$ in $(\Gamma_{G} \setminus \{ 1 \})^c$. So assume $g_1 \in \langle w \rangle$. Since $\langle w \rangle$, contain a unique cyclic subgroup of order $d_1$, we have $g_2 \not\in \langle w \rangle$, and then $g_2$ is adjacent to $w$ in $(\Gamma_{G} \setminus \{ 1 \})^c$.
 
\bigskip

\noindent {\bf Case II.} $d_1 \mid e_m$

\bigskip

\noindent We use the same notations as in Case I. Here, if $w \not\in \langle g_1 \rangle$, we are done. So assume that $w \in \langle g_1 \rangle$. Now, the number of cyclic classes that are involved in $\gamma_{d_1}$ is $\geq 3$ and only one of them is connected to $\gamma_{d_2}$ inside ${\mathcal{H}}$. So, we can reconstruct ${\mathcal{H}}$ by changing its initial vertex $w$ by $w^{\prime}$, where $w^{\prime}$ belongs to another cyclic class of order $d_1$. Then, $w^{\prime} \not\in \langle g_1 \rangle$ and we are done.  

\bigskip

\noindent Once this step is done, the remaining vertices in Step 2 can be joined exactly the way they were in Step 1. We denote by ${\mathcal{H}}^{\prime}$ the extended Hamilton path. 

\bigskip

\noindent {\bf Step 3.}

\bigskip

\noindent Now we consider the vertices in $g \in \cup_{C \in {\mathfrak{S}}} C$ so that there is a unique cyclic class consisting of elements of order $|g|$ inside $\cup_{C \in {\mathfrak{S}}} C$ and ${\mathfrak{m}}_G(|g|) \geq 2$. We denote this cyclic class by $D$. Then, $|g| = d_j$ for some $1 \leq j \leq k$ (see Step 1). Then these vertices are adjacent to the vertices in $\gamma_j$ in $(\Gamma_{G} \setminus \{ 1 \})^c$. Since the number of cyclic classes consisting of elements of order $d_j$ is $\geq 3$, we have $|V(\gamma_{d_j})| \geq 3|D|$. Then we remove the old edges in the Hamilton path ${\mathcal{H}}$ and attach the vertices of $D$ as given in the following diagram:

\begin{figure}[htp] \centering{
\includegraphics[scale=0.42]{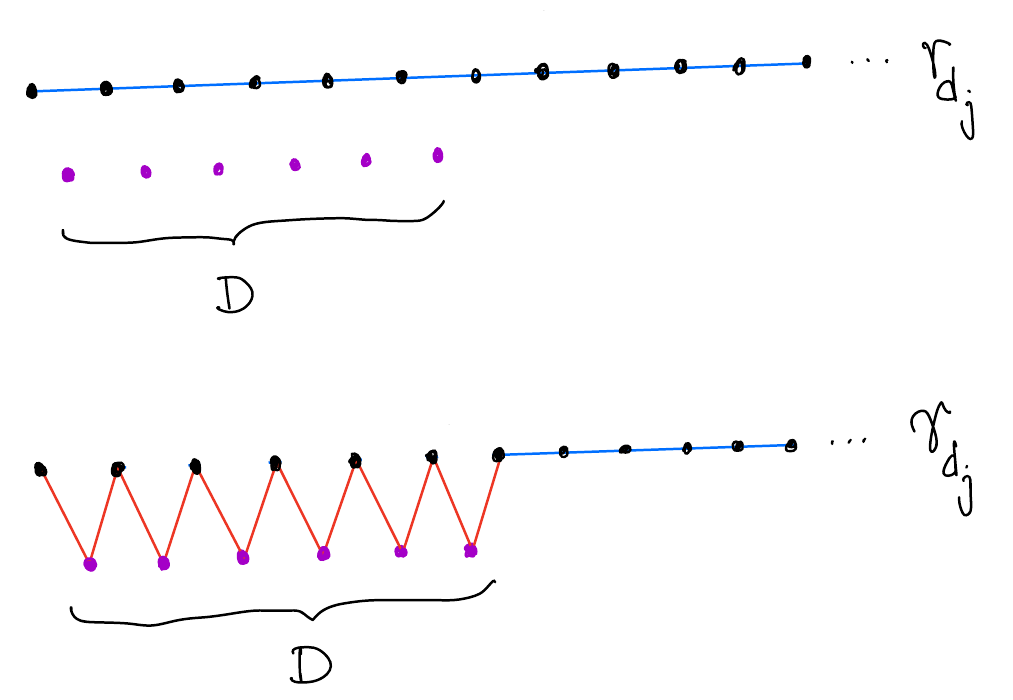}}
\caption{}
\end{figure} 

\vspace*{.25in}

\noindent We denote the modified Hamilton path by ${\mathcal{H}}^{\prime\prime}$. 

\bigskip

\noindent {\bf Step 4.} 

\bigskip

\noindent The remaining elements $g \in \cup_{C \in {\mathfrak{S}}} C$ satisfies ${\mathfrak{m}}_G(|g|) = 1$. Let $w$ be a vertex of ${\mathcal{H}}$ (i.e., from $G \setminus \cup_{C \in {\mathfrak{S}}} C$), so that either $g \in \langle w \rangle$, or $w \in \langle g \rangle$. Then, we have $\langle w \rangle \cap U(1;G) \neq 1$, a contradiction. So, $g$ is adjacent to $w$ in $(\Gamma_{G} \setminus \{ 1 \})^c$. We can now attach these vertices $g$ to the vertices of ${\mathcal{H}}$, which is incident to at most one vertex considered in Step 3 (see the figure below):  

\begin{figure}[htp] \centering{
\includegraphics[scale=0.52]{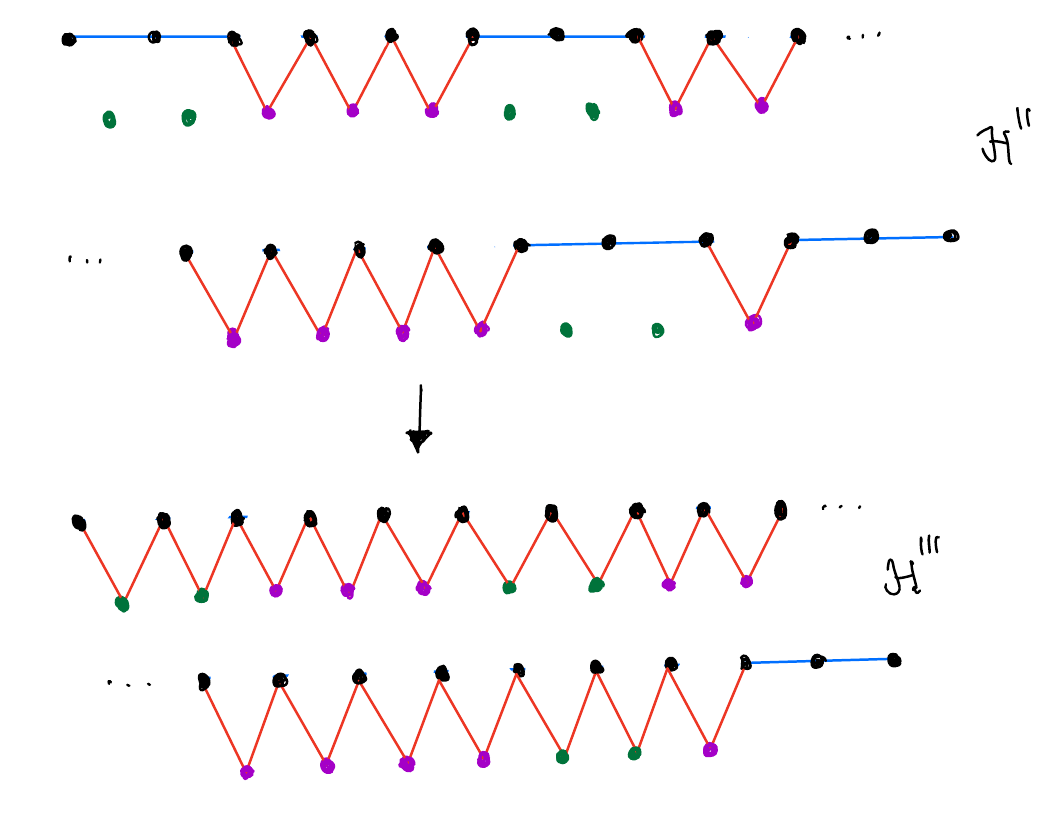}}
\caption{}
\end{figure} 

\bigskip

\noindent We need to prove that there are enough such available vertices in $G \setminus \cup_{C \in {\mathfrak{S}}} C$. We note that  
\begin{eqnarray*}
\lvert \big\{ g \in G \setminus \{ 1 \} ~:~ {\mathfrak{m}}_G(|g|) = 1 \big\} \rvert & \leq & \lvert U(1;G) \setminus \{ 1 \} \rvert ~~<~~ \lvert U(1;G) \rvert \\
	& \leq & \lvert \cup_{C \in {\mathfrak{S}}} C \rvert \leq \lvert G \setminus \cup_{C \in {\mathfrak{S}}} C \rvert
\end{eqnarray*}
from the hypothesis. In particular, if $\langle \cup_{C \in {\mathfrak{S}}} C \rangle$ is a proper subgroup of $G$, then we have 
\[
\lvert \cup_{C \in {\mathfrak{S}}} C \rvert \leq \lvert \langle \cup_{C \in {\mathfrak{S}}} C \rangle \rvert \leq \lvert G \setminus \langle \cup_{C \in {\mathfrak{S}}} C \rangle \rvert \leq \lvert G \setminus \cup_{C \in {\mathfrak{S}}} C \rvert,
\]
since $[G : \langle \cup_{C \in {\mathfrak{S}}} C \rangle] \geq 2$. This proves that the modified Hamilton path in $(\Gamma_{G} \setminus \{ 1 \})^c$ in covers all vertices of $(\Gamma_{G} \setminus \{ 1 \})^c$, and the proof is complete. \QED

\bigskip

\section{A family of groups $G$ for which $\lambda(G) > |G|$}
\label{final-example}

\bigskip

\noindent Let $p$ be an odd prime $\geq 5$, and consider the direct product $E(p) := C_p \times {\mathbb{D}}_6 = \langle x, \rho, \sigma \rangle$, where $C_p = \langle x \rangle$ denote the cyclic group of order $p$, and ${\mathbb{D}}_6 = \langle \rho, \sigma \rangle$ is the dihedral group of order $6$ with $\rho^3  = 1 = \sigma^2$. Then $\Delta(E(p)) = \{ 1, 2, 3, p, 2p, 3p \}$, and ${\mathrm{Fit}}(E(p)) = U(1; E(p)) = \langle x, \rho \rangle \cong C_{3p}$. Here we note that $E(p)$ is not nilpotent, since it does not have a normal Sylow-$2$-subgroup, and it is full uni-core. 

\bigskip

\noindent The cyclic subgroups $C \leq E(p)$ with $C \cap U(1; E(p)) \neq 1$ are given in the following table: 
\begin{center}
\begin{tabular}{ |c|c| } 
 \hline
$\lvert C \rvert$ & $C$ \\ \hline
$3$ & $\langle \rho \rangle$ \\ \hline
$p$ & $\langle x \rangle$ \\ \hline
$3p$ & $\langle x, \rho \rangle$ \\ \hline 
$2p$ & $\langle x, \sigma \rangle, \langle x, \sigma\rho \rangle, \langle x, \sigma\rho^2 \rangle$ \\ \hline 
\end{tabular}
\end{center}

\noindent Note that the only cyclic subgroups that are disjoint from $U(1;G)$ are the three subgroups of order $2$ contained in ${\mathbb{D}}_6$. However, these subgroups are already contained in some cyclic subgroup of order $2p$. Hence, $\Delta(E(p)) \setminus \Delta_0(E(p)) = \{ 2 \}$, even though $\cup_{C \in {\mathfrak{S}}} C = E(p)$. Hence, the condition (ii) in Theorem \ref{main-thm-classification} is false. The final conclusion is given as:

\bigskip

\begin{theorem}\label{excep-example} 
(a) If $p=5$, then $\lambda(E(5)) = \lvert E(5) \rvert$, and

\smallskip

\noindent (b) for all prime $p \geq 7$, we have $\lambda(E(p)) > \lvert E(p) \rvert$. 
\end{theorem}

\bigskip

\noindent {\bf Proof.} We set $G := E(p)$.

\smallskip

\noindent (a) We first look at the list of cyclic classes $F$ in $G \setminus \{ 1 \}$ given in the following table:
\begin{center}
\begin{tabular}{ |c|c|c| } 
 \hline
$\lvert g \rvert, g \in F$ & $F$ & No. of classes \\ \hline
$2$ & $\{ \sigma \}$, $\{ \sigma\rho \}$, $\{ \sigma\rho^2 \}$ & $3$ \\ \hline
$3$ & $\{ \rho, \rho^2 \}$ & $1$ \\ \hline
$5$ & $\{ x, x^2, x^3, x^4 \}$ & $1$ \\ \hline
$10$ & $\{ x\sigma, x^2 \sigma, x^3 \sigma, x^4 \sigma \}$, $\{ x\sigma\rho, x^2 \sigma\rho, x^3 \sigma\rho, x^4 \sigma\rho \}$, $\{ x\sigma\rho^2, x^2 \sigma\rho^2, x^3 \sigma\rho^2, x^4 \sigma\rho^2 \}$ & $3$ \\ \hline 
$15$ & $\{ x^i \rho^j ~:~ i = 1,2,3,4; j = 1,2 \}$ & $1$ \\ \hline 
\end{tabular}
\end{center}  

\noindent Now let $F_1, F_2, F_3$ denote the cyclic classes consisting of elements of order $10$. Using the construction in section \ref{main-thm-classification}, we have a Hamiltonian path $\gamma(F_1, F_2, F_3) = (w_1, w_2, \dotsc, w_{12})$ in the subgraph of $G \setminus \{ 1 \}$ induced by $F_1 \cup F_2 \cup F_3$. If $J = \{ v_1, \dotsc, v_8 \}$ denote the cyclic class consisting of elements of order $15$, then we modify $\gamma(F_1, F_2, F_3)$ to the Hamiltonian path ${\mathcal{H}} = (w_1, v_1, w_2, v_2, \dotsc, w_8, v_8, w_9, w_{10}, w_{11}, w_{12})$ in the subgraph of $G \setminus \{ 1 \}$ induced by $F_1 \cup F_2 \cup F_3 \cup J$.  Then, ${\mathcal{H}}$ can be extended to the Hamiltonian path ${\mathcal{H}}^{\prime} = {\mathcal{H}} (\rho, x, \rho^2, x^2, \sigma, x^3, \sigma\rho, x^4, \sigma\rho^2)$ which cover all vertices of $G \setminus \{ 1 \}$.

\bigskip

\noindent (b) Using Theorem \ref{mfw-main-lemma}, it is enough to prove that $(\Gamma_{G} \setminus \{ 1 \})^c$ does not admit a Hamiltonian path. If possible, assume that such a path ${\mathcal{H}}$ exists. Let $g \in G$ be an element of order $p$. Then, $g$ is adjacent to all elements of $\langle x \rangle$ in the power graph $\Gamma_G$. We claim that $g$ is adjacent to any element of order $2p$ or $3p$ in $\Gamma_G$. Let $h \in G$ be an element of order $2p$ (resp. $3p$). Then $h^2$ (resp. $h^3$) has order $p$ and hence there exists $i \in {\mathbb{N}}$ with ${\mathrm{g.c.d.}}(i,p) = 1$ so that $(h^2)^i = g$ (resp. $(h^3)^i = g$). This implies that $g$ and $h$ are adjacent in $\Gamma_G$. Thus, in the Hamiltonian path ${\mathcal{H}} \subseteq (\Gamma_{G} \setminus \{ 1 \})^c$, $g$ can be only be adjacent to the elements of order $2$ or $3$, which are given by $\{ \rho, \rho^2, \sigma, \sigma\rho, \sigma\rho^2 \}$. This set has size $5 < 6 \leq p-1$, since $p \geq 7$. Then, by the Pigeonhole principle, there exists an element of order $p$ which must be adjacent to an element of order $p, 2p$ or $3p$ in the complement graph $(\Gamma_{G} \setminus \{ 1 \})^c$,  a contradiction. \QED      

\bigskip
 
\bibliographystyle{plain} 
\bibliography{Uniwidth-univ-lambda}

@article {akc,
    AUTHOR = {Abawajy, Jemal and Kelarev, Andrei and Chowdhury, Morshed},
     TITLE = {Power graphs: a survey},
   JOURNAL = {Electron. J. Graph Theory Appl. (EJGTA)},
  FJOURNAL = {Electronic Journal of Graph Theory and Applications. EJGTA},
    VOLUME = {1},
      YEAR = {2013},
    NUMBER = {2},
     PAGES = {125--147},
   MRCLASS = {05C25 (05C40 05C45 05C78)},
  MRNUMBER = {3145411},
       DOI = {10.5614/ejgta.2013.1.2.6},
       URL = {https://doi.org/10.5614/ejgta.2013.1.2.6},
}

@article {ksc,
    AUTHOR = {Kumar, Ajay and Selvaganesh, Lavanya and Cameron, Peter J. and
              Tamizh Chelvam, T.},
     TITLE = {Recent developments on the power graph of finite groups---a
              survey},
   JOURNAL = {AKCE Int. J. Graphs Comb.},
  FJOURNAL = {AKCE International Journal of Graphs and Combinatorics},
    VOLUME = {18},
      YEAR = {2021},
    NUMBER = {2},
     PAGES = {65--94},
      ISSN = {0972-8600},
   MRCLASS = {05C25 (05C40 05C45 05C78)},
  MRNUMBER = {4310374},
       DOI = {10.1080/09728600.2021.1953359},
       URL = {https://doi.org/10.1080/09728600.2021.1953359},
}

@article {cgs,
    AUTHOR = {Chakrabarty, Ivy and Ghosh, Shamik and Sen, M. K.},
     TITLE = {Undirected power graphs of semigroups},
   JOURNAL = {Semigroup Forum},
  FJOURNAL = {Semigroup Forum},
    VOLUME = {78},
      YEAR = {2009},
    NUMBER = {3},
     PAGES = {410--426},
      ISSN = {0037-1912},
   MRCLASS = {20M99 (05C25 05C40 05C45)},
  MRNUMBER = {2511776},
MRREVIEWER = {V\'{a}clav Koubek},
       DOI = {10.1007/s00233-008-9132-y},
       URL = {https://doi.org/10.1007/s00233-008-9132-y},
}

@incollection {kq,
    AUTHOR = {Kelarev, A. V. and Quinn, S. J.},
     TITLE = {A combinatorial property and power graphs of groups},
 BOOKTITLE = {Contributions to general algebra, 12 ({V}ienna, 1999)},
     PAGES = {229--235},
 PUBLISHER = {Heyn, Klagenfurt},
      YEAR = {2000},
   MRCLASS = {05C25 (05C20 20F99)},
  MRNUMBER = {1777663},
MRREVIEWER = {M. E. Watkins},
}

@incollection {falc,
    AUTHOR = {Fern\'{a}ndez-Alcober, Gustavo A.},
     TITLE = {An introduction to finite {$p$}-groups: regular {$p$}-groups and groups of maximal class},
      NOTE = {16th School of Algebra, Part I (Portuguese) (Bras\'{\i}lia, 2000)},
   JOURNAL = {Mat. Contemp.},
  FJOURNAL = {Matem\'{a}tica Contempor\^{a}nea},
    VOLUME = {20},
      YEAR = {2001},
     PAGES = {155--226},
      ISSN = {0103-9059},
   MRCLASS = {20D15},
  MRNUMBER = {1868828},
MRREVIEWER = {Eamonn A. O'Brien},
}

@article {mfw,
    AUTHOR = {Ma, Xuanlong and Feng, Min and Wang, Kaishun},
     TITLE = {Lambda number of the power graph of a finite group},
   JOURNAL = {J. Algebraic Combin.},
  FJOURNAL = {Journal of Algebraic Combinatorics. An International Journal},
    VOLUME = {53},
      YEAR = {2021},
    NUMBER = {3},
     PAGES = {743--754},
      ISSN = {0925-9899},
   MRCLASS = {05C78 (05C25)},
  MRNUMBER = {4258067},
       DOI = {10.1007/s10801-020-00940-9},
       URL = {https://doi.org/10.1007/s10801-020-00940-9},
}

@article {mil,
    AUTHOR = {Miller, G. A.},
     TITLE = {An {E}xtension of {S}ylow's {T}heorem},
   JOURNAL = {Proc. London Math. Soc. (2)},
  FJOURNAL = {Proceedings of the London Mathematical Society. Second Series},
    VOLUME = {2},
      YEAR = {1905},
     PAGES = {142--143},
      ISSN = {0024-6115},
   MRCLASS = {DML},
  MRNUMBER = {1577264},
       DOI = {10.1112/plms/s2-2.1.142},
       URL = {https://doi.org/10.1112/plms/s2-2.1.142},
}

@article {kula,
    AUTHOR = {Kulakoff, A.},
     TITLE = {\"{U}ber die {A}nzahl der eigentlichen {U}ntergruppen und der {E}lemente von gegebener {O}rdnung in {$p$}-{G}ruppen},
   JOURNAL = {Math. Ann.},
  FJOURNAL = {Mathematische Annalen},
    VOLUME = {104},
      YEAR = {1931},
    NUMBER = {1},
     PAGES = {778--793},
      ISSN = {0025-5831},
   MRCLASS = {DML},
  MRNUMBER = {1512698},
       DOI = {10.1007/BF01457969},
       URL = {https://doi.org/10.1007/BF01457969},
}

@book {isa,
    AUTHOR = {Isaacs, I. Martin},
     TITLE = {Character theory of finite groups},
      NOTE = {Corrected reprint of the 1976 original [Academic Press, New
              York; MR0460423]},
 PUBLISHER = {AMS Chelsea Publishing, Providence, RI},
      YEAR = {2006},
     PAGES = {xii+310},
      ISBN = {978-0-8218-4229-4; 0-8218-4229-3},
   MRCLASS = {20-02 (20C05 20C15 20C20 20C25)},
  MRNUMBER = {2270898},
       DOI = {10.1090/chel/359},
       URL = {https://doi.org/10.1090/chel/359},
}

@article {ber,
    AUTHOR = {Berkovi\v{c}, Ja. G.},
     TITLE = {{$p$}-groups of finite order},
   JOURNAL = {Sibirsk. Mat. \v{Z}.},
  FJOURNAL = {Akademija Nauk SSSR. Sibirskoe Otdelenie. Sibirski\u{\i} Matemati\v{c}eski\u{\i} \v{Z}urnal},
    VOLUME = {9},
      YEAR = {1968},
     PAGES = {1284--1306},
      ISSN = {0037-4474},
   MRCLASS = {20.40},
  MRNUMBER = {0241534},
MRREVIEWER = {O. H. Kegel},
}

@book {ber2,
    AUTHOR = {Berkovich, Yakov},
     TITLE = {Groups of prime power order. {V}ol. 1},
    SERIES = {De Gruyter Expositions in Mathematics},
    VOLUME = {46},
      NOTE = {With a foreword by Zvonimir Janko},
 PUBLISHER = {Walter de Gruyter GmbH \& Co. KG, Berlin},
      YEAR = {2008},
     PAGES = {xx+512},
      ISBN = {978-3-11-020418-6},
   MRCLASS = {20D15 (20-02)},
  MRNUMBER = {2464640},
MRREVIEWER = {Eamonn A. O'Brien},
       DOI = {10.1515/9783110208238.512},
       URL = {https://doi.org/10.1515/9783110208238.512},
}

@article {cgh,
    AUTHOR = {Cameron, Peter J. and Ghosh, Shamik},
     TITLE = {The power graph of a finite group},
   JOURNAL = {Discrete Math.},
  FJOURNAL = {Discrete Mathematics},
    VOLUME = {311},
      YEAR = {2011},
    NUMBER = {13},
     PAGES = {1220--1222},
      ISSN = {0012-365X},
   MRCLASS = {05C25},
  MRNUMBER = {2793235},
       DOI = {10.1016/j.disc.2010.02.011},
       URL = {https://doi.org/10.1016/j.disc.2010.02.011},
}

@article {cam,
    AUTHOR = {Cameron, Peter J.},
     TITLE = {The power graph of a finite group, {II}},
   JOURNAL = {J. Group Theory},
  FJOURNAL = {Journal of Group Theory},
    VOLUME = {13},
      YEAR = {2010},
    NUMBER = {6},
     PAGES = {779--783},
      ISSN = {1433-5883},
   MRCLASS = {05C25 (20D60)},
  MRNUMBER = {2736156},
MRREVIEWER = {Mohammad A. Iranmanesh},
       DOI = {10.1515/JGT.2010.023},
       URL = {https://doi.org/10.1515/JGT.2010.023},
}

@article {ms2,
    AUTHOR = {Sarkar, Siddhartha and Mishra, Mayank},
     TITLE = {The lambda number of the power graph of a finite {$p$}-group},
   JOURNAL = {J. Algebraic Combin.},
  FJOURNAL = {Journal of Algebraic Combinatorics. An International Journal},
    VOLUME = {57},
      YEAR = {2023},
    NUMBER = {1},
     PAGES = {101--110},
      ISSN = {0925-9899,1572-9192},
   MRCLASS = {05C25 (05C78 20D15)},
  MRNUMBER = {4544260},
       DOI = {10.1007/s10801-022-01158-7},
       URL = {https://doi.org/10.1007/s10801-022-01158-7},
}

\end{document}